\documentclass[12pt]{article}
\usepackage{graphicx} 

\usepackage{epsfig}
\usepackage{amsmath}
\usepackage{amssymb}
\usepackage{amsthm}
\usepackage{latexsym}
\usepackage{cite}

\newcommand{\comment}[1]{}
\newtheorem{theorem}{Theorem}        
\newtheorem{remark}{Remark}
\newtheorem{corollary}{Corollary}

\begin{document}
 
\title{\LARGE
{\bf The Percolation Transition in the Zero-Temperature Domany Model}
}

\author{
{\bf Federico Camia}
\thanks{Forschungsinstitut f\"ur Mathematik, ETH, 8092 Z\"urich, Switzerland}\,
\thanks{Present address: EURANDOM, P.O. Box 513, 5600 MB Eindhoven, The Netherlands.}\,
\thanks{E-mail: camia@eurandom.tue.nl}\\
\and
{\bf Charles M.~Newman}
\thanks{Courant Inst.~of Mathematical Sciences, New York University, New York, NY 10012, USA}\,
\thanks{Research partially supported by the
U.S. NSF under grant DMS-01-04278.}\,
\thanks{E-mail: newman@courant.nyu.edu}
}

\date{}

\maketitle

\begin{abstract}
We analyze a deterministic cellular automaton 
$\sigma^{\cdot} = (\sigma^n : n \geq 0)$ corresponding to the
zero-tem\-per\-a\-ture case of Domany's stochastic Ising ferromagnet on
the hexagonal lattice $\mathbb H$.
The state space ${\cal S}_{\mathbb H} = \{ -1, +1 \}^{\mathbb H}$ consists
of assignments of $-1$ or $+1$ to each site of $\mathbb H$ and the initial
state $\sigma^0 = \{ \sigma_x^0 \}_{x \in {\mathbb H}}$ is chosen
randomly with $P(\sigma_x^0 = +1) = p \in [0,1]$.
The sites of $\mathbb H$ are partitioned in two sets $\cal A$ and $\cal B$
so that all the neighbors of a site $x$ in $\cal A$ belong to $\cal B$ and
vice versa, and the discrete time dynamics is such that the $\sigma^{\cdot}_x$'s
with $x \in {\cal A}$ (respectively, $\cal B$) are updated simultaneously
at odd (resp., even) times, making $\sigma^{\cdot}_x$ agree with the
majority of its three neighbors.

In \cite{cns} it was proved that there is a percolation transition at
$p = 1/2$ in the percolation models defined by $\sigma^n$, for all times
$n \in [1, \infty]$.
In this paper, we study the nature of that transition and prove that
the critical exponents $\beta$, $\nu$ and $\eta$ of the dependent
percolation models defined by $\sigma^n, n \in [1, \infty]$, have the
same values as for standard two-dimensional independent 
site percolation (on the triangular lattice).
\end{abstract}



\section{Introduction}
The deterministic cellular automaton corresponding to the zero-temperature
case of Domany's stochastic Ising ferromagnet on the hexagonal lattice $\mathbb H$
\cite{domany} can be considered as a simplified version of a continuous
time Markov process where an independent (rate $1$) Poisson clock is assigned
to each site $x \in {\mathbb H}$, and the spin $\sigma_x$ at site $x$ is
updated when the corresponding clock rings.
The rule for updating the spin is to flip it if and only if it disagrees with
two or three (a majority) of its neighbors.
This model has been studied both rigorously and numerically in~\cite{hn}; the
numerical results about two critical 
exponents obtained there strongly suggest that the
dependent percolation model defined by the limiting spin configuration
$\sigma^{\infty}$ is in the same universality class as ordinary independent
percolation.

In the zero-temperature Domany model studied here, the rule for updating the
spins is unchanged, but the timing is different.
This model has been previously used in numerical simulations~\cite{nienhuis},
and its dynamical as well as percolation properties have been studied in
\cite{cns}, \cite{cns1} and \cite{cns2}.
In~\cite{cns}, the existence of a percolation transition for the dependent
percolation models associated with the spin configurations $\sigma^n$ at time
$n$, for all $n \in [1, \infty]$, is proved; in~\cite{cns1}, it is shown that
for any $n \geq 1$, the
crossing probabilities converge to Cardy's formula~\cite{cardy} when the
lattice spacing $\delta$ is sent to zero (the \emph{continuum scaling
limit}), as in the case of independent percolation
(at least on the triangular lattice $\mathbb T$~\cite{smirnov}); in~\cite{cns2},
the continuum scaling limit is analyzed in terms of cluster boundaries and is
shown to be the same as for ordinary independent critical percolation
on $\mathbb T$.
This last result strongly suggests that also the critical exponents
defined at the critical point should be the same as for independent
percolation. 
Indeed, due to the above mentioned result, one can use properties of the
Stochastic Loewner Evolution~\cite{schramm} (which describes the scaling
limit of ordinary critical percolation cluster interfaces) to compute
certain critical exponents in the continuous model corresponding to the scaling limit.
But unfortunately, the connection between critical exponents in the continuous
and discrete models is not straightforward; more work and further results
on the discrete model would be required to relate the discrete critical
exponents with the continuous ones (see~\cite{sw}).

Similar models on different lattices have been studied in various papers;
see, for example,~\cite{cdn, fss, gns, nns, ns1, ns2, ns3}  for models on
${\mathbb Z}^d$ and~\cite{Howard} for a model on the homogeneous tree of
degree three.
Such models are also discussed  extensively in the physics literature,
usually on ${\mathbb Z}^d$ (see, for example,~\cite{domany} and~\cite{lms}).
Their interest is tied to the fact that they can be obtained as zero-temperature
limits of stochastic Ising models, a special class of Markov processes
whose transition probabilities/rates are chosen so that the Gibbs
measures (for some Hamiltonian) at temperature $T$ are invariant for the
Markov process.
In systems where there are multiple (infinite-volume) Gibbs measures for
$T$ below some critical $T_c$, a subject of considerable interest is the
$t \to \infty$ behavior of the spin configuration $\sigma^t$ (with
temperature $T = T_1 < T_c$) when the initial state is chosen from the
(unique) Gibbs measure at $T = T_2 > T_c$.
Studying the limiting case where $T_1 = 0$ and $T_2 = \infty$ is the standard
choice in much of the statistical physics literature (see, e.g., \cite{bray}).

In this paper, we study the percolation properties of (say) $+1$ spins in
the dependent percolation models, generated by the zero-temperature Domany
cellular automaton, corresponding to $\sigma^n$, with $n \in [1, \infty]$,
when the initial state is chosen randomly with $P(\sigma_x^0 = +1) = p \in [0,1]$.
In~\cite{cns}, it was shown that there is a percolation transition with
$p_c = 1/2$ for all values of $n \in [1, \infty]$.
We remark that $p_c$ is \emph{not} $1/2$ for $n=0$, since the critical
probability for independent site percolation on the hexagonal lattice is
strictly larger than $1/2$, therefore the system is \emph{driven}
to criticality by the dynamics, and this after just one time step.
This is in contrast to what happens in the case of the continuous time
model studied in~\cite{hn}, where it is believed that criticality is
achieved again for $p=1/2$, but only at time $t = \infty$.
Nonetheless, the nature of the percolation transition,
in terms of critical exponents and/or scaling limits, is presumed
to be the same in the different models.

\section{Definition of the Model and Preliminary Results} \label{model}

In this section, we give a more detailed description of the model
and present, for completeness, results that were proved in~\cite{cns}
which motivate (and will be used in) the next section, where the
main results of this paper are presented.

\begin{figure}[!ht]
\begin{center}
\includegraphics[width=8cm]{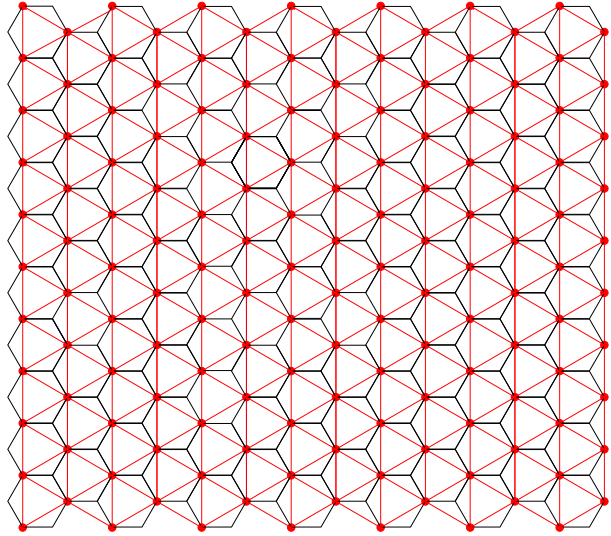}
\caption{A star-triangle transformation relating triangular and hexagonal lattice.}
\label{star-tri}
\label{int}
\end{center}
\end{figure}

Consider the homogeneous ferromagnet on the hexagonal lattice
${\mathbb H}$ (embedded in ${\mathbb R}^2$ so that the elementary
cells are regular hexagons with side length $1$ -- see, for example,
Figure~\ref{star-tri}) with states denoted by 
$\sigma = \{ \sigma_x \}_{x \in {\mathbb H}}, \, \sigma_x = \pm 1$, 
and with (formal) Hamiltonian
\begin{equation} \label{Hamiltonian}
{\cal H}(\sigma) = - \sum_{ \langle x,y \rangle } \sigma_x \sigma_y ,
\end{equation}
where $\sum_{ \langle x,y \rangle }$ denotes the sum over all pairs of 
neighbor sites, each pair counted once.
We write ${\cal N}_{\mathbb H}(x)$ for the set of 
three neighbors of $x$, and indicate with
\begin{equation} \label{variation}
\Delta_x  {\cal H} (\sigma) =
2 \sum_{y \in {\cal N}_{\mathbb H}(x)} \sigma_x \sigma_y
\end{equation}
the change in the Hamiltonian when the \emph{spin} $\sigma_x$
at site $x$ is flipped (i.e., changes sign).

Partition the sites of the hexagonal lattice ${\mathbb H}$ into two subsets,
${\cal A}$ and ${\cal B}$, in such a way that all three neighbors of a site
$x$ in ${\cal A}$ (resp., ${\cal B}$) are in ${\cal B}$ (resp., ${\cal A}$).
By joining two sites of ${\cal A}$ whenever they are next-nearest neighbors
in the hexagonal lattice (two steps away from each other),
we get a triangular lattice (the same with ${\cal B}$ -- see Figure~\ref{star-tri}).
The synchronous dynamics that we consider here is such that all the sites in
${\cal A}$ (resp., ${\cal B}$) are updated simultaneously.

We now define the cellular automaton $\sigma^n, \, n \in {\mathbb N}$,
with  state space ${\cal S}_{\mathbb H} = \{ -1,+1 \}^{\mathbb H}$,
which is the zero temperature limit of a model of Domany~\cite{domany},
as follows:
\begin{itemize}
\item The initial state $\sigma^0$ is chosen 
from a Bernoulli product measure $P_p$, with $P_p(\sigma^0_0 = +1) = p$.
\item At odd times $n = 1, 3, \dots$, 
the spins in ${\cal A}$ are updated according to the
following rule: $\sigma_x, \, x \in {\cal A}$, 
is flipped if and only if $\Delta_x  {\cal H} (\sigma)<0$.
\item At even times $n = 2, 4, \dots$, 
the spins in ${\cal B}$ are updated according to the
same rule as for those in ${\cal A}$.
\end{itemize}

We denote by $\sigma^{\infty}$ the limiting state of the
cellular automaton $\sigma^n$.
$\sigma^{\infty} = \lim_{n \to \infty} \sigma^n$ exists with
probability one~\cite{nns} and, like $\sigma^n$
for $1 \leq n < \infty$, defines a dependent percolation model
on $\mathbb H$.

The following observations are useful in understanding the behavior
of the model and will help in the proof of 
our main result, Theorem~\ref{mainthm-critexp},
which is presented in the next section of the paper.
\begin{itemize}
\item The values of the spins in $\cal A$ at time $0$ are
irrelevant since after the first update those values are
uniquely determined by the values of the spins in $\cal B$. 
\item Once the initial spin configuration in $\cal B$ is
chosen, the dynamics is completely deterministic.
\item A spin can no longer flip once it belongs to either a loop or
a ``barbell'' of constant sign, where a barbell consists of two
disjoint loops connected by a path.
\end{itemize}

There is an alternative but equivalent way of describing the
discrete time dynamics as a deterministic cellular automaton
(with random initial state) on the triangular lattice $\mathbb T$
(corresponding to the set $\cal B$ of sites).
Given some site $\bar{x} \in {\mathbb T}$, group its six ${\mathbb T}$-neighbors 
$y_i^{\bar x}$ in three disjoint pairs $\{ y_1^{\bar{x}}, y_2^{\bar{x}} \}$,
$\{ y_3^{\bar{x}}, y_4^{\bar{x}} \}$, $\{ y_5^{\bar{x}}, y_6^{\bar{x}} \}$, so that
$y_1^{\bar{x}}$ and $y_2^{\bar{x}}$ are ${\mathbb T}$-neighbors,
and so on for the other two pairs.
Translate this construction to all sites $x \in {\mathbb T}$, thus producing
three pairs of sites  $\{ y_1^x, y_2^x \}$, $\{ y_3^x, y_4^x \}$,
$\{ y_5^x, y_6^x \}$ associated to each site $x \in {\mathbb T}$.
(Note that this construction does not need to specify how
${\mathbb T}$ is embedded in ${\mathbb R}^2$.)
Site $x$ is updated at times $m = 1, 2, \ldots$ according to the following rule:
the spin at site $x$ is changed from $\sigma_x$ to $- \sigma_x$ if and only if 
at least two of its pairs of neighbors have all four sites with the same sign
$- \sigma_x$.

The models on the hexagonal and on the triangular lattice are related
through a \emph{star-triangle transformation}
(see Figure~\ref{star-tri} and, for example, p.~335~of~\cite{grimmett}).
More precisely, the dynamics on the triangular lattice ${\mathbb T}$ is
equivalent to the zero-temperature Domany dynamics on the hexagonal lattice
${\mathbb H}$ when restricted to the sublattice ${\cal B}$ for even times $n=2m$.

To see this, start with ${\mathbb T}$ and construct an
hexagonal lattice ${\mathbb H}'$ by means of a star-triangle transformation
such that a site is added at the center of each of the triangles
$(x, y_1^x, y_2^x), (x, y_3^x, y_4^x)$, and $(x, y_5^x, y_6^x)$.
${\mathbb H}'$ may be partitioned into two triangular sublattices ${\cal A}'$
and ${\cal B}'$ with ${\cal B}' = {\mathbb T}$.
One can now see that the dynamics on ${\mathbb T}$ for $m=1,2,\ldots$
and the zero-temperature Domany dynamics on ${\mathbb H}'$ restricted to
${\cal B}'$ for even times $n=2m$ are the same.

An immediate consequence of this equivalence between the two cellular automata
is that the two families of percolation models that they produce are also
equivalent in an obvious way through a star-triangle transformation.
To be more precise, the percolation models defined on ${\mathbb T}$
for times $m=1,2,\ldots$ are the same as those defined on ${\cal B}$ by
the zero-temperature Domany model for even times $n=2m$.

We now present (without proof) the results of~\cite{cns} on the zero-temperature
Domany model.
Theorem~\ref{syn_exp} says that the convergence to the limiting state is
exponentially fast, while Theorem~\ref{mean_c_s} identifies the critical point
of the percolation transition, $p_c = 1/2$.

\begin{theorem} \label{syn_exp}
Let $P_{\cal A} (n)$ denote the probability that a deterministic site in
${\cal A}$ flips after time $n$ and similarly for $P_{\cal B} (n)$.
Then, for any $p$, there is a constant $c \in (0, \infty)$ such that
\begin{equation} \label{syn_expon}
P_{\cal A} (n), P_{\cal B} (n) \leq e^{-cn}.
\end{equation}
\end{theorem}

\begin{theorem} \label{mean_c_s}
If $p > 1/2$ (resp., $< 1/2$), there is percolation of $+1$ (resp., $-1$)
spins in $\sigma^n$ for any $n \in [1, \infty]$ for almost every $\sigma^0$. \\

\noindent If $p = 1/2$ \emph{:}
\begin{enumerate}
\item For $n \in [0, \infty]$, there is no percolation in $\sigma^n$ of either $+1$ or
$-1$ spins, for almost every $\sigma^0$.
\item The mean cluster size is infinite in $\sigma^n$ for all $n \in [1, \infty]$
(but finite for $n=0$).
\end{enumerate}
\end{theorem}

\section{Critical Exponents} \label{critexp}

We will consider three percolation critical exponents, namely the
exponents $\beta$ (related to the percolation probability), $\nu$
(related to the correlation length) and $\eta$ (related to the connectivity function).
The existence of these exponents has been proved, and their predicted values
confirmed rigorously, in recent papers~\cite{lsw,sw}, for the case of independent site
percolation on the triangular lattice.
Such exponents are believed to be universal for independent percolation in
the sense that their value should depend only on the number of dimensions and not on
the structure of the lattice or on the nature of the percolation model 
(e.g., whether it is site or bond percolation);
that type of universality has not yet been proved.

Consider an independent percolation model with distribution $P_p$ on a two-dimensional
lattice $\mathbb L$ such that $0<p_c<1$.
Let $C_x$ be the open cluster containing site $x$ and $|C_x|$ be its cardinality,
then $\theta(p) = \theta_x(p) = P_p(|C_x| = \infty)$ is the \emph{percolation probability}.
Arguments from theoretical physics suggest that $\theta(p)$ behaves roughly like
$(p-p_c)^{\beta}$ as $p$ approaches $p_c$ from above.

It is also believed that the
\emph{connectivity function}
\begin{equation} \label{connect-function}
\tau_p(x,y) = P_p(x \text{ and } y \text{ belong to the same cluster})
\end{equation}
behaves, for the Euclidean 
length $||x-y||$ large, like $||x-y||^{-\eta}$ if $p=p_c$, and like
$\exp{(-||x-y||/\xi(p))}$ if $0 < p < p_c$,
for some $\xi(p)$ satisfying $\xi(p) \to \infty$ as $p \uparrow p_c$.
The \emph{correlation length} $\xi(p)$ is defined by
\begin{equation} \label{correlation-length}
\xi(p)^{-1} = \lim_{||x-y|| \to \infty} \left\{ - \frac{1}{||x-y||} \log \tau_p(x,y) \right\}.
\end{equation}
$\xi(p)$ is expected to behave like $(p_c - p)^{-\nu}$ as $p \uparrow p_c$.

It is not clear how strong one may expect such asymptotic relations to be
(for more details about critical exponents and \emph{scaling theory} in
percolation, see~\cite{grimmett} and references therein).
For this reason, the logarithmic relation is usually employed.
That means that the previous conjectures are usually stated in the following form:
\begin{eqnarray}
\lim_{p \downarrow p_c} \frac{\log \theta(p)}{\log (p-p_c)} = \beta, \\
\lim_{||x-y|| \to \infty} \frac{\log \tau_{p_c}(x,y)}{\log ||x-y||} = - \eta, \\
\lim_{p \uparrow p_c} \frac{\log \xi(p)}{\log (p_c - p)} = - \nu.
\end{eqnarray}

In the rest of this section and in the next one, $\theta_x(p)$, $\tau_p(x,y)$
and $\xi(p)$ will indicate the percolation probability, connectivity function
and correlation length for independent site percolation on the triangular lattice
$\mathbb T$ (identified with $\cal B$).
We will denote by $\theta_x(p,n)$, $\tau_{p,n}(x,y)$ and $\xi(p,n)$ the corresponding
quantities for the percolation models on the hexagonal lattice $\mathbb H$ at time
$n = 0, 1, 2, \ldots$ .
The main theorem of this paper is the following.

\begin{theorem} \label{mainthm-critexp}
There exist constants $0 < c_1, c_2, c_3, c_4 < \infty$ such that, $\forall n \in [1, \infty]$,
and $x, y \in {\mathbb H}$ and suitably chosen $x',y' \in {\cal B}$ with
$||x-x'||, ||y-y'|| \leq 1$,
\begin{eqnarray}
c_1 \, \theta_{x'}(p) \leq \theta_x(p,n) \leq c_2 \, \theta_{x'}(p), \,\,\,
\text{ for } p \in (1/2,1], \label{theta} \\
p^{c_3} \, \tau_p(x',y') \leq \tau_{p,n}(x,y) \leq p^{-c_4} \, \tau_p(x',y'), \,\,\,
\text{ for } p \in (0, 1/2], \label{tau} \\
\xi(p,n) = \xi(p), \,\,\, \text{ for } p \in (0, 1/2]. \label{xi}
\end{eqnarray}
\end{theorem}

The next corollary is an immediate consequence of Theorem~\ref{mainthm-critexp}
and its main application; it says that the dependent percolation models defined
by $\sigma^n$, with $n \in [1, \infty]$, are in the universality class of ordinary
independent percolation.

\begin{corollary} \label{cor-critexp}
The critical exponents  $\beta$, $\eta$ and $\nu$ exist for the dependent
percolation models on $\mathbb H$ defined by $\sigma^n$, with $n \in [1, \infty]$,
and have the same numerical values as for independent site percolation on $\mathbb T$.
\end{corollary}

\begin{remark} As already mentioned, the existence of the exponents $\beta$, $\nu$
and $\eta$ for independent site percolation on $\mathbb T$ has been recently proved,
and their predicted values confirmed rigorously~\cite{lsw, sw}.
\end{remark}

\section{Proofs}

Before we can prove the main results of this paper, we need some notation.
We will denote by $P_{p,n}$ the distribution of $\sigma^n$ with initial density
$p$ of plus spins (i.e., with initial distribution $P_{p,0} = P_p$).

For a site $x \in {\mathbb H}$, we will denote by ${\cal N}_{\mathbb H}(x)$ the set
of its three neighbors in $\mathbb H$.
For a site $x \in {\cal B}$, we will denote by ${\cal N}_{\cal B}(x)$ the set
of its six neighbors in $\cal B$ endowed with the graph structure of a triangular
lattice (as explained in Section~\ref{model} -- see also Figure~\ref{star-tri}).

We will call an \emph{$\mathbb H$-path} a path on the hexagonal lattice
$\mathbb H$ and a \emph{$\cal B$-path} a path on the triangular lattice $\cal B$.
A path whose spins are all plus (resp., minus) will be called a plus (resp., minus) path.
Similarly, we will call an \emph{$\mathbb H$-loop} a (simple) loop on the hexagonal lattice
$\mathbb H$ and \emph{$\cal B$-loop} a (simple) loop on the triangular lattice $\cal B$.
A loop whose spins are all plus (resp., minus) will be called a plus (resp., minus) loop.
Notice that constant-sign $\mathbb H$-loops, doubly-infinite $\mathbb H$-paths and
``barbells'' (two loops connected by a path) are stable for the dynamics, in the
sense that, once formed, their spins will never flip again.

For two subsets $C$ and $D$ of $\mathbb H$, we indicate with
$\{C \stackrel{\mathbb H}{\longleftrightarrow} D\}$
the event that some site in $C$ is connected to some site in $D$ by a \emph{plus}
$\mathbb H$-path, with $\{C \stackrel{\mathbb H}{\longleftrightarrow} \infty\}$
the event that some site in $C$ belongs to an infinite \emph{plus}
$\mathbb H$-path.
For two subsets $C$ and $D$ of $\cal B$, endowed with the graph structure of
a triangular lattice, we indicate with $\{C \stackrel{\cal B}{\longleftrightarrow} D\}$
the event that some site in $C$ is connected to some site in $D$ by a \emph{plus} 
$\cal B$-path,
with $\{C \stackrel{\cal B}{\longleftrightarrow} \infty\}$ 
the event that some site in $C$ belongs to an infinite \emph{plus} $\cal B$-path.

\bigskip

\noindent {\it Proof of Theorem~\ref{mainthm-critexp}.}
Let us first assume that $x$ and $y$ belong to $\cal B$; we then take $x' = x, y' = y$
in~(\ref{theta}) -- (\ref{tau}).
The lower bound for $\theta_x(p,n)$ in Eq.~(\ref{theta}) comes from the following
bound
\begin{eqnarray}
\theta_x(p,n) & = & P_{p,n}(x \stackrel{\mathbb H}{\longleftrightarrow} \infty) \\
& \geq & P_{p,0}(\{ x \stackrel{\cal B}{\longleftrightarrow} \infty \} \cap
\{ x \text{ belongs to a plus $\mathbb H$-loop} \}) \\
 & \geq &  P_{p,0}(x \stackrel{\cal B}{\longleftrightarrow} \infty) \,
P_{p,0}(x \text{ belongs to a plus $\mathbb H$-loop}) \\
 & \geq & p^6 \, \theta_x(p) \\
 & > & \frac{1}{2^6} \, \theta_x(p),
\end{eqnarray}
where we have used the fact that at time $1$ (when the sites in $\cal A$
are updated for the first time) the dynamics transforms any constant sign
$\cal B$-path into a constant sign $\mathbb H$-path,
the FKG inequality and the fact that the events
$\{ 0 \stackrel{\cal B}{\longleftrightarrow} \infty \}$ and
$\{ 0 \text{ belongs to a plus $\mathbb H$-loop} \}$ are increasing (see, for example,
\cite{grimmett}), and the fact that the smallest $\mathbb H$-loop contains $6$ sites.

The lower bound for $\tau_{p,n}(x,y)$ in Eq.~(\ref{tau}) is obtained in a
similar way.
Mimicking the proof of the lower bound in Eq.~(\ref{theta}), we have immediately
\begin{equation}
\tau_{p,n}(x,y) = P_{p,n}(x \stackrel{\mathbb H}{\longleftrightarrow} y) \geq p^{12}
\, \tau_p(x,y).
\end{equation}

For the upper bound of Eq.~(\ref{theta}) we rely on the following observation.
If no site in ${\cal N}_{\cal B}(x)$ belongs to an infinite plus $\cal B$-path
at time $0$, then, by the self-matching property of the triangular lattice, site
$x$ must be surrounded by a minus $\cal B$-loop, which will produce a stable
$\mathbb H$-loop at time $1$.
Therefore, site $x$ will not belong to an infinite plus $\mathbb H$-path at any later time.
Thus,
\begin{equation} \label{upper-theta}
\theta_x(p,n) \leq P_{p,0}({\cal N}_{\cal B}(x) \stackrel{\cal B}{\longleftrightarrow} \infty).
\end{equation}

Since the event $\{ x \stackrel{\cal B}{\longleftrightarrow} \infty \}$ can be written as
$\{ \sigma^0_x = +1 \} \cap \{ {{\cal N}_{\cal B}}(x) \stackrel{\cal B}{\longleftrightarrow} \infty \}$,
using the FKG inequality we have
\begin{equation}
P_{p,0}(x \stackrel{\cal B}{\longleftrightarrow} \infty) \geq
p \, P_{p,0}({\cal N}_{\cal B}(x) \stackrel{\cal B}{\longleftrightarrow} \infty).
\end{equation}
From this and Eq.~(\ref{upper-theta}), we get
\begin{equation}
\theta_x(p,n) \leq p^{-1} \, \theta_x(p) \leq 2 \, \theta_x(p),
\end{equation}
as required.

%

For the upper bound of Eq.~(\ref{tau}), we first note that for bounded $||x-y||$,
the inequality is trivial by choosing $c_4$ big enough so that the right-hand side
of~(\ref{tau}) exceeds $1$.
Next, for $||x-y||$ large enough, we notice that unless
$\{ {\cal N}_{\cal B}(x) \stackrel{\cal B}{\longleftrightarrow} {\cal N}_{\cal B}(y) \}$
at time $0$, $x$ and $y$ must be separated by a minus $\cal B$-loop surrounding
one of them or by a doubly-infinite minus $\cal B$-path, and therefore it cannot be the
case that $\{x \stackrel{\mathbb H}{\longleftrightarrow} y\}$ at any later time because
at time $1$ a stable minus $\mathbb H$-loop or doubly-infinite $\mathbb H$-path will be formed.
This yields, for $||x-y||$ large enough,
\begin{equation} \label{upper-tau}
\tau_{p,n}(x,y) \leq
P_{p,0}({\cal N}_{\cal B}(x) \stackrel{\cal B}{\longleftrightarrow} {\cal N}_{\cal B}(y)).
\end{equation}

Since the event $\{x \stackrel{\cal B}{\longleftrightarrow} y\}$ can be written as
$\{ \sigma^0_x = \sigma^0_y = +1 \} \cap
\{ {\cal N}_{\cal B}(x) \stackrel{\cal B}{\longleftrightarrow} {\cal N}_{\cal B}(y) \}$,
using the FKG inequality we have
\begin{equation}
P_{p,0}(x \stackrel{\cal B}{\longleftrightarrow} y) \geq
p^2 \, P_{p,0}({\cal N}_{\cal B}(x) \stackrel{\cal B}{\longleftrightarrow} {\cal N}_{\cal B}(y)).
\end{equation}
From this and Eq.~(\ref{upper-tau}), we get
\begin{equation}
\tau_{p,n}(x,y) \leq p^{-2} \, \tau_p(x,y),
\end{equation}
as required.

If $x$ and $y$ belong to $\cal A$, the proof is analogous, but one has to consider
slightly different events.
In this case we take $x' \in {\cal N}_{\mathbb H}(x)$ and then we have, for the lower bound
in Eq.~(\ref{theta}),
\begin{eqnarray}
\theta_x(p,n) & = & P_{p,n}(x \stackrel{\mathbb H}{\longleftrightarrow} \infty) \\
& \geq & P_{p,0}(\{ x' \stackrel{\cal B}{\longleftrightarrow} \infty \}
\cap \{ x \text{ belongs to a plus $\mathbb H$-loop} \}) \\
 & \geq &  P_{p,0}(x' \stackrel{\cal B}{\longleftrightarrow} \infty) \,
P_{p,0}(x \text{ belongs to a plus $\mathbb H$-loop}) \\
 & \geq & p^6 \, \theta_{x'}(p) \\
 & > & \frac{1}{2^6} \, \theta_{x'}(p),
\end{eqnarray}
where $||x-x'||=1$ because $x'$ is an $\mathbb H$-neighbor of $x$.

Analogously, for the lower bound in Eq.~(\ref{tau}), we get
\begin{equation}
\tau_{p,n}(x,y) = P_{p,n}(x \stackrel{\mathbb H}{\longleftrightarrow} y) \geq p^{12}
\, \tau_p(x',y'),
\end{equation}
with $x' \in {\cal N}_{\mathbb H}(x)$ and $y' \in {\cal N}_{\mathbb H}(y)$.

For the upper bound of Eq.~(\ref{theta}), we notice that, if no site in
$\cup_{z \in {\cal N}_{\mathbb H}(x)} {\cal N}_{\cal B}(z)$ belongs to an
infinite plus $\cal B$-path at time $0$, then, by the self-matching
property of the triangular lattice, site $x$ must be surrounded by a minus
$\cal B$-loop that will produce a stable $\mathbb H$-loop at time $1$.
Therefore, site $x$ will not belong to an infinite plus $\mathbb H$-path at any later time.
Thus,
\begin{equation}
\theta_x(p,n) \leq P_{p,0}(\{ \cup_{z \in {\cal N}_{\mathbb H}(x)} {\cal N}_{\cal B}(z)
\stackrel{\cal B}{\longleftrightarrow} \infty).
\end{equation}

Since for $x' \in {\cal N}_{\mathbb H}(x)$,
$\{ x' \stackrel{\cal B}{\longleftrightarrow} \infty \}
\supset \{ \sigma^0_z = +1,  \, \forall z \in {\cal N}_{\mathbb H}(x) \} \cap
\{ \cup_{z \in {\cal N}_{\mathbb H}(x)} {\cal N}_{\cal B}(z) \stackrel{\cal B}{\longleftrightarrow} \infty \}$,
using the FKG inequality we have
\begin{equation}
P_{p,0}(x' \stackrel{\cal B}{\longleftrightarrow} \infty) \geq
p^3 \, P_{p,0}(\cup_{z \in {\cal N}_{\mathbb H}(x)} {\cal N}_{\cal B}(z) \stackrel{\cal B}{\longleftrightarrow} \infty),
\end{equation}
which yields
\begin{equation}
\theta_x(p,n) \leq p^{-3} \, \theta_{x'}(p) \leq 8 \, \theta_{x'}(p),
\end{equation}
with $||x-x'||=1$.

The upper bound of Eq.~(\ref{tau}), still trivial for bounded $||x-y||$,
follows again from a similar observation for $||x-y||$ large enough: unless
$\{ \cup_{z \in {\cal N}_{\mathbb H}(x)} {\cal N}_{\cal B}(z) \stackrel{\cal B}{\longleftrightarrow}
\cup_{z' \in {\cal N}_{\mathbb H}(y)} {\cal N}_{\cal B}(z') \}$
at time $0$, $x$ and $y$ must be separated by a minus $\cal B$-loop surrounding
one of them or by a doubly-infinite minus $\cal B$-path.
At time $1$, a stable minus $\mathbb H$-loop or doubly-infinite $\mathbb H$-path
will be formed, making the event
$\{x \stackrel{\mathbb H}{\longleftrightarrow} y\}$ impossible at any later time.
This yields, for $||x-y||$ large enough,
\begin{equation}
\tau_{p,n}(x,y) \leq
P_{p,0}(\cup_{z \in {\cal N}_{\mathbb H}(x)} {\cal N}_{\cal B}(z)
\stackrel{\cal B}{\longleftrightarrow} \cup_{z' \in {\cal N}_{\mathbb H}(y)} {\cal N}_{\cal B}(z')).
\end{equation}

Since for $x' \in {\cal N}_{\mathbb H}(x)$ and $y' \in {\cal N}_{\mathbb H}(y)$,
$\{ x' \stackrel{\cal B}{\longleftrightarrow} y' \}
\supset \{ \sigma^0_z = +1, \, \forall z \in {\cal N}_{\mathbb H}(x) \} \cap
\{ \sigma^0_{z'} = +1, \, \forall z' \in {\cal N}_{\mathbb H}(y) \}
\cap \{ \cup_{z \in {\cal N}_{\mathbb H}(x)} {\cal N}_{\cal B}(z) \stackrel{\cal B}{\longleftrightarrow}
\cup_{z' \in {\cal N}_{\mathbb H}(y)} {\cal N}_{\cal B}(z') \}$,
using the FKG inequality we have
\begin{equation}
P_{p,0}(x' \stackrel{\cal B}{\longleftrightarrow} y') \geq
p^6 \, P_{p,0}(\cup_{z \in {\cal N}_{\mathbb H}(x)} {\cal N}_{\cal B}(z)
\stackrel{\cal B}{\longleftrightarrow} \cup_{z' \in {\cal N}_{\mathbb H}(y)} {\cal N}_{\cal B}(z')),
\end{equation}
which yields
\begin{equation}
\tau_{p,n}(x,y) \leq p^{-6} \, \tau_p(x',y'),
\end{equation}
with $x' \in {\cal N}_{\mathbb H}(x)$ and $y' \in {\cal N}_{\mathbb H}(y)$.

The proof of Eq.~(\ref{tau}) in the remaining case (namely, $x \in {\cal A}, \, y \in {\cal B}$
or vice versa) should now be clear.

Eq.~(\ref{xi}) is an immediate consequence of Eq.~(\ref{tau}) and the definition of $\xi(p)$;
it is enough to observe that
\begin{equation}
\lim_{||x - y|| \to \infty} \left\{ - \frac{1}{||x' - y'||} \left[\log \tau_p(x',y') + c_3 \log p \right] \right\}
= \xi(p)^{-1}
\end{equation}
and
\begin{equation}
\lim_{||x - y|| \to \infty} \left\{ - \frac{1}{||x' - y'||} \left[\log \tau_p(x',y') - c_4 \log p \right] \right\}
= \xi(p)^{-1}. \,\,\,\,\, \fbox{}
\end{equation}

\bigskip

\noindent {\it Proof of Corollary~\ref{cor-critexp}.}
It follows from Eqs.~(\ref{theta})~and~(\ref{tau}) that, for $p \in (1/2,1]$ and $||x-y||$ large enough,
\begin{eqnarray}
- \frac{\log c_1 + \log \theta_{x'}(p)}{\log (p - 1/2)} \leq - \frac{\log \theta_x(p,n)}{\log (p - 1/2)}
\leq - \frac{\log c_2 + \log \theta_{x'}(p)}{\log (p - 1/2)}, \\
\frac{\log \tau_{1/2}(x',y') + c_3 \log \frac{1}{2}}{\log ||x-y||} \leq \frac{\log \tau_{1/2, n}(x,y)}{\log ||x-y||}
\leq \frac{\log \tau_{1/2}(x',y') - c_4 \log \frac{1}{2}}{\log ||x-y||}.
\end{eqnarray}
Using these two equations, together with Eq.~(\ref{xi}) and the definitions of the critical
exponents, and taking the appropriate limits gives the desired results. \fbox{}


%

\end{document}